\documentclass[a4paper,leqno,11pt,twoside]{amsart}
\usepackage{amsmath,amsfonts,amssymb,amsthm,amscd}
\usepackage[utf8]{inputenc}
\usepackage[english]{babel}
\usepackage[colorlinks=true,citecolor=blue, urlcolor=blue, linkcolor=blue,pagebackref]{hyperref}
\usepackage[top=1.1in,bottom=1.1in,left=1.in,right=1.in]{geometry} 
\usepackage{graphicx}
\usepackage{wrapfig}
\usepackage{paralist}
\usepackage{tabto}
\usepackage{standalone}
\usepackage{xfrac}
\usepackage{float}
\usepackage{tikz-cd}
\usepackage{mathtools}
\usepackage{cleveref}

\usepackage{pgf,tikz}
\usetikzlibrary{arrows,chains,
	decorations.pathreplacing,
	shapes,%
}

\usepackage[all]{xy} 

\usepackage{enumerate}
\usepackage{xcolor}
\usepackage{aliascnt}

\theoremstyle{plain}
\newtheorem{thm}{Theorem}[section]

\newtheorem{thmIntr}{Theorem}

\newaliascnt{propIntr}{thmIntr}

\newaliascnt{corIntr}{thmIntr}

\newaliascnt{QU}{thm}
\newtheorem{QU}[QU]{Question}
\aliascntresetthe{QU}

\newaliascnt{lem}{thm}
\newtheorem{lem}[lem]{Lemma}
\aliascntresetthe{lem}

\newaliascnt{cor}{thm}
\newtheorem{cor}[cor]{Corollary}
\aliascntresetthe{cor}

\newaliascnt{prop}{thm}
\newtheorem{prop}[prop]{Proposition}
\aliascntresetthe{prop}

\theoremstyle{definition}

\newaliascnt{rem}{thm}
\newtheorem{rem}[rem]{Remark}
\aliascntresetthe{rem}

\newaliascnt{defn}{thm}

\aliascntresetthe{defn}

\newaliascnt{ex}{thm}

\aliascntresetthe{ex}

\numberwithin{equation}{section}

\def\bP{\ensuremath{\mathbb{P}}}

\def\bZ{\ensuremath{\mathbb{Z}}}
\def\bC{\ensuremath{\mathbb{C}}}

\def\cI{\ensuremath{\mathcal{I}}}

\def\cO{\ensuremath{\mathcal{O}}}

\def\cW{\ensuremath{\mathcal{W}}}

\newcommand{\GalCl}[0]{W} 
\newcommand{\degrm}[0]{d} 
\newcommand{\RatMap}[0]{\gamma} 
\newcommand{\Subgr}[0]{G} 
\newcommand{\EllCurve}[0]{C} 
\newcommand{\ShortKer}[0]{K} 

\definecolor{applegreen}{rgb}{0.55, 0.71, 0.0}

\title[An explicit class of Lagrangian surfaces]{An explicit class of Lagrangian surfaces}
\author[P.~Grossi and F.~Moretti]{Paolo Grossi and Federico Moretti}
\thanks{Paolo Grossi is a member of GNSAGA (INdAM)}

\address{Dipartimento di Matematica ”F. Casorati”,
Universit\`{a} di Pavia, Via Adolfo Ferrata 5, 27100 Pavia, Italy}
\email{paolo.grossi01@universitadipavia.it} 

\address{Stony Brook University,  100 Nicolls Road, Stony Brook, NY 11794}
\email{federico.moretti@stonybrook.edu}

\begin{document}

\begin{abstract}
We construct a family of general type surfaces with $q=4$, $p_g=6$ and $K^2=24$. These surfaces enjoy some interesting properties: they are Lagrangian in their Albanese variety and their canonical map is $2:1$ onto a  degree $12$ surface in $\mathbb P^5$ with $44$ even nodes. 
\end{abstract}
	

\maketitle

\setcounter{tocdepth}{1}

\section{Introduction}

\quad 
Let $(T,\eta)$ be a $2n-$dimensional variety equipped with a non degenerate $2-$form $\eta$. An $n-$dimensional subvariety $X\subset T$ is called Lagrangian if $\eta_{|_X}=0$.
In this paper, we consider $W$, the Galois closure of a degree $3$ rational map $\varphi \colon A \dashrightarrow \mathbb P^2$ for a very general abelian surface $(A,L)$ of type $(1,6)$. We obtain an explicit Lagrangian surface in $A\times A$ with respect to the invariant $2-$form under the action of the Galois group. This construction has already been carried out in \cite{BPS} for a $(1,2)$ abelian surface. We compute the invariants of $W$ and study the geometry of the canonical map and of the Albanese map. Moreover, we obtain that the rational map $\varphi \colon A\dashrightarrow \mathbb P^2$ is equivariant with respect to a group of translations $K\simeq \mathbb Z_3\oplus \mathbb Z_3$. Denoting $X=W/K$, the central result of the paper is the following. 
	
\begin{thmIntr}
    The surface $X$ has invariants
    \[K_{X}^2 = 24, \quad c_2(X) = 12, \quad \chi(\cO_{X}) = 3, \quad q(X) = 4, \quad p_g = 6.\]
    Moreover the following hold:
    \begin{itemize}
        \item the canonical map $X\to \mathbb P^5$ is $2:1$ onto a degree $12$ surface with $44$ even nodes;
        \item the surface $X$ is Lagrangian into its Albanese variety, the Albanese map $\alpha: X\to (A \times A)/K$ is a local immersion, there exists a unique couple of distinct points $(p,q)\in X\times X$ such that $\alpha(p)=\alpha(q)$.
    \end{itemize}
    \label{ThmA}
\end{thmIntr}
\quad The construction is governed by a unique canonical genus $6$, hyperelliptic curve in $|L|$ (which turns out to be nodal for general $A$). The following result may be of independent interest.	

\begin{thmIntr}
    Given  a very general $(1,6)$ abelian surface $(A,L)$ there is a unique hyperelliptic curve of geometric genus $6$ in $|L|$ (up to the action of $K(L_A)$).
    \label{ThmB}
    \label{introhyp}\end{thmIntr}

 \quad The surfaces we construct in \autoref{ThmA} are closely related to those constructed by Schoen in \cite{schoen2007} and studied by various authors, see for example \cite{cilimlroll}, \cite{ritorousarti}, \cite{beavilleshoen}, \cite{campana}. Schoen surfaces are of general type with invariants $q=4,p_g=5,K^2=16$ (so they are Lagrangian in their Albanese varieties). The canonical map of a Schoen surface is a $2:1$ cover of a $(2,4)$ complete intersection in $\mathbb P^4$ with $40$ even nodes. Beauville considered in \cite{beavilleshoen} some covers of Schoen surfaces obtaining surfaces with invariants $q=4,p_g=7,K^2=32$ whose canonical map is $2:1$ onto a $(2,2,2,2)$ complete intersection in $\mathbb P^6$ with $48$ even nodes. Our surface completes the tale in $\mathbb P^5$ (quoting \cite{beavilleshoen}): it fits in this list of examples with canonical map of degree $2$ onto a canonical low degree surface of general type  with an high number of nodes. Notice also that in all these examples $K^2=8\chi.$  Let us remark that by \cite{beauvillemappacanonica} the canonical image of a surface of general type is either rational or a canonically embedded surface, the second case being quite exceptional if the canonical map is of positive degree onto its image (for a list of different examples where this happens see \cite{cilipartove}).  In the recent survey \cite{survey} several interesting questions are raised concerning the degree of the canonical map of general type surfaces.

\quad This list of surfaces is interesting also by a topological point of view. The fact that the map \[
\lambda: \bigwedge ^2H^0(\Omega^1_X)\to H^0(\omega_X)
\]
has a non-simple tensor in the kernel has strong implications on the structure of the fundamental group: it is non-commutative and one can ask whether it is nihilpotent, see \cite{cam}. In \cite{campana} the property of being Lagrangian (or more generally \emph{coisotropic}) in an ambient variety is systematically investigated. In our case (inspired by \cite{BPS}) the Lagrangian structure comes from the action of the Galois group $S_3$ on $(W,A\times A)$, whose quotient is (birationally) $(\mathbb P^2,\mathrm{Kum}_2(A))$, where $\mathrm{Kum}_2A$ denotes the Kummer four-fold associated to $A$. Since $\mathrm{Kum}_2(A)$ carries a $2-$form  and $\mathbb P^2$ is Lagrangian, the pullback of this form to $A\times A$ vanishes on $W$. We are also able to show that the rank of $\lambda $ is exactly $5$. In particular $\lambda$ does not have any simple tensor in the kernel, hence by Castelnuovo-de Franchis theorem $X$ does not admit any map to a curve of genus $\ge 2$. Other examples of Lagrangian surfaces in abelian four-folds were constructed in \cite{bogotschi}.

\quad The original motivation of Schoen was to prove the Hodge conjecture for some abelian four-folds.  Indeed, for a general Schoen surface $S$, the cohomology class of $S$ in $H^{2,2}(\mathrm{Alb}(S))$ is not a complete intersection class. This allowed Schoen to prove the Hodge conjecture for a general abelian four-fold of type $\mathrm{Alb}(S)$ for a Schoen surface $S$.

\quad Some natural questions remain open for the family of surfaces we construct. \begin{QU}
    What is the dimension of the moduli space of such surfaces?
\end{QU}
\quad We constructed such a surface for a very general abelian surface of type $(1,6)$, so the dimension of the moduli space is at least $3$. It is natural to wonder whether these surfaces can be deformed such that their Albanese is no longer isogenous to a product. A natural way to tackle this problem is to use the approach of Schoen in \cite{schoen2007}, i.e. applying semiregularity techniques \cite{bloch}, \cite{semiregolarita2} (for more details see \autoref{hodge}). As mentioned before, this may prove the Hodge conjecture for a class of abelian four-folds. Another natural approach would be to compute the expected dimension of $(2,2,3)$ complete intersections in $\mathbb P^5$ with $44$ even nodes (supposing our surfaces are complete intersections), as done in \cite{cilimlroll} to compute alternatively the dimension of the moduli space of Schoen surfaces. Unfortunately, in our case, the expected dimension is $2$. This motivates our following. 
\begin{QU}
    Are these surfaces $(2,2,3)$ complete intersections in $\mathbb P^5$? Can one find their equations?
\end{QU}
\quad This program has been carried out in \cite{ritorousarti} for Schoen surfaces. By \cite[Theorem 3.1]{konno} (see also \cite[Theorem 4.2]{konno}) the canonical image is a $(2,2,3)$ complete intersection unless it is contained in a rational normal scroll. We suspect that both situations may occur (the latter being a specialization of the former).

\quad As mentioned before, these surfaces are interesting from a topological perspective. 
\begin{QU}
    What is the topology of these surfaces? Is the nihilpotent tower of their fundamental group finite? What is the universal cover?
\end{QU}
\quad Understanding the structure of the fundamental group is likely to be very challenging.

\quad \textbf{Structure of the paper.} \Cref{Sectio: 2} is devoted to study the geometry of the degree \(3\) rational map $\varphi:A\dashrightarrow \mathbb P^2$. In \Cref{Section: 3} we analyze the case of such rational maps on quotients of products of elliptic curves. We also show that the hyperelliptic curve of \autoref{ThmB} is nodal of genus $6$. In \Cref{Section: 4} we show that the ramification of the resolution of $\varphi$ is a cover of a genus \(6\) smooth hyperelliptic curve and we compute its arithmetic genus and self-intersection. In \Cref{Section: 5}, this is the main tool to compute the invariants of $X$ the Galois closure of $\varphi$. The last sections are devoted to the study of the geometry of this surface: \Cref{sec:AlbaneseMap} concerns its Albanese map, in \Cref{Section: 7} we show that it is not fibered over curves of positive genus and in \Cref{sec:CanonicalImage} we study its canonical map.

\quad \textbf{Acknowledgements.} We would like to thank Pietro Pirola for suggesting this project and for numerous fruitful discussions. We thank Frederick Campana for interesting discussions on the fundamental group of these surfaces. The authors thank Rita Pardini for some interesting comments and questions. We thank the referee for several suggestions.

\section{The degree \texorpdfstring{\(3\)}{3} cover}
\label{Sectio: 2}

The construction of the degree $3$ cover $A \dashrightarrow \mathbb P^2$ has been sketched by the second author in \cite{mor}. The goal of this section is to investigate and understand explicitly this rational map. From now on we shall assume that $(A,L_A)$ is a $(1,6)$ abelian surface. Every such surface admits a description as $A=B/\langle b \rangle $, 
where $(B,L_B)$ is a $(1,3)$ abelian surface and $b$ is a non trivial $2-$torsion point. Now let us denote by $\pi:B \to A$ the degree $2$ isogeny, we get that $F=\pi_*L_B$ is a rank $2$ vector bundle with $c_2(F)=3$. \\ 
First of all, we report the result of \cite{mor} (to which we address the reader for a complete proof). We will suppose that the line bundle $L_B$ does not admit a subsheaf  of the form $\mathcal O_B(nC)$ where $C$ is an elliptic curve and $n\ge 2$ and that all divisors in $|L_B|$ do not contain an elliptic curve invariant by translating by $b$. This is a sufficient condition for $L_B$ to be globally generated (see \cite[Lemma 10.1.2]{BirkenhakeLange1992complex}) and for $F$ to be globally generated in codimension $1$, which is the condition required in \cite{mor} to be able to construct a degree $3$ morphism (see Section 1, Section 3 and Lemma 3.1 of \cite{mor}).

\begin{thm}
   \label{degree3cover}
Let \((B,L_B)\), \(b \in B\), \(A=B/ \langle b \rangle\) and \(F=\pi_*L_B\) be as above. Then, the exact sequence 
    \begin{equation}
    \begin{tikzcd}
    0 \arrow{r} & F^\vee \arrow{r} & H^0(F)^\vee \otimes \cO_A \arrow{r} & L_A,
    \end{tikzcd}
    \label{ev}
    \end{equation}
induces a degree \(3\) (rational) map \(\varphi_{H^0(F)^\vee}:A \dashrightarrow \mathbb PH^0(F)=\mathbb P^2\), the fiber above the class of \(s\in H^0(F)\) being \(Z(s)\). 
 
\end{thm}
\begin{proof} We outline the proof, for more details we refer the reader to \cite[Section 3]{mor}.
 As soon as $L_B$ does not admit a  subsheaf of the form $\mathcal O_B(nC)$, where $C$ is an elliptic curve and $n\ge 2$, and no divisors in $|L_B|$ contain an elliptic curve invariant by $t_b$, the vector bundle $F=\pi_*L_B$ is globally generated in codimension $1$ (see \cite[Lemma 3.1]{mor}). This, in turn, implies that the image of the evaluation map $\mathrm{ev}:H^0(F)\otimes\mathcal O_A\to F$ has the same first Chern class as $F$ and that $\mathrm{ker}(\mathrm{ev})=L_A^\vee$ (it is a line bundle with first Chern class $-c_1(F))$. Dualizing one gets an exact sequence \[
 \begin{tikzcd}
     0\arrow{r}& F^\vee \arrow{r}{\mathrm{ev}^*}  & H^0(F)^\vee \otimes \mathcal O_A \arrow{r} & L_A,
 \end{tikzcd}
 \]
 with the last map surjective in codimension $1$ (in the same locus where $\mathrm{ev}$ is surjective). From the above sequence one gets an inclusion $H^0(F)^\vee\subset H^0(L_A)$. The linear system $H^0(F)^{\vee}$ induces the map $\varphi_{H^0(F)^\vee}:A\dashrightarrow \mathbb PH^0(F)$. The base locus coincides with the locus where $\mathrm{ev}$ is not surjective. Then, it is easy to check that after removing the base locus of the map, the sequence (\ref{ev}) is the pullback of the twisted Euler sequence \[
\begin{tikzcd}
    0\arrow{r} & (T_{\mathbb P H^0(F)}(-1))^\vee \arrow{r} & H^0(F)^\vee \otimes \mathcal O_{\mathbb PH^0(F)}\arrow{r} & \mathcal O(1) \arrow{r} & 0.
\end{tikzcd}
\]
Now any section $s\in H^0(T_{\mathbb PH^0(F)}(-1))=H^0(F)$ vanishes in exactly one point $[s]$ and the fiber above the class of $[s]\in \mathbb PH^0(F)$ is $Z(\varphi_{H^0(F)^\vee}^*(s))$ which is a section of $F$ (slightly abusing notation, we called this section $s$ in the statement of the theorem).

\end{proof}
	
To simplify notation we will call:
	
\begin{itemize}
    \item \(\varphi \coloneqq \varphi_{H^0(F)^\vee}\) defined in \autoref{degree3cover};
    \item \(\cI\) the ideal of the base locus of \(\varphi\);
    \item \(\tilde{\varphi} \colon \tilde{A} \to \bP^2\) the resolution of singularities of \(\varphi\);
    \item \((-1)_A \colon A \to A\) (respectively \((-1)_B \colon B \to B\)) the morphism mapping an element of \(A\) (respectively \(B\)) to its opposite;
    \item \(K(L_B)\) the kernel of the natural isogeny \(B \to \mbox{Pic}^0(B)\) mapping \(x \in B\) to \(L_B^{-1} \otimes t_x^{*} L_B\) (\(t_x\) being the translation by \(x\))\footnote{Starting from \autoref{sec:ramification}, in order to further simplify notation, we will refer to it as \(K\).}.
\end{itemize}
	
Investigating \(\varphi\), we get some symmetry results.
	
\begin{thm}
In the condition of \autoref{degree3cover}, the following hold:
    \begin{itemize}
        \item \(\varphi\) is equivariant with respect to the action of \(K(L_B)\); 
	\item \(\varphi\) is equivariant with respect to \(-1\);
        \item the base locus of \(\varphi\) is \(a+K(L_B)\) where \(a\in A[2]\) is a \(2-\)torsion point.
    \end{itemize}
(By abuse of notation, we identify \(K(L_B)\) and \(\pi(K(L_B))\)).
\label{equivariance}
\end{thm}
	
\begin{proof}
    To prove the equivariance with respect to the action of \(K(L_B)\), for $x\in K(L_B)$ let us denote with $t_x,t_{\pi(x)}$ the translations by $x$ and $\pi(x)$ respectively. Since the following diagram
        \[
        \begin{tikzcd}
        B\arrow{r}{t_x} \arrow{d}{\pi}& B \arrow{d}{\pi} \\
        A \arrow{r}{t_{\pi(x)}} &A
        \end{tikzcd}
        \]
    is cartesian, we get an isomorphism $t_{\pi(x)}^*F=t_{\pi(x)}^*\pi_*L_B\simeq \pi_*t_{x}^*L_B\simeq \pi_*L_B$, from which the equivariance follows.\\
    The equivariance with respect to \(-1\) follows by a similar argument on the commutative diagram 
        \[
        \begin{tikzcd}
        B\arrow{r}{(-1)_B} \arrow{d}{\pi}& B \arrow{d}{\pi} \\
        A \arrow{r}{(-1)_A} &A.
        \end{tikzcd}
        \]
    once we notice that it is always possible to choose the origin of \(B\) so that \((-1)^{*}_BL_B \simeq L_B\).\\
    Finally, to prove the third point of the statement, we notice that the number of base points is finite and, in fact, they are at most \(9\). Indeed, \(c_2(L \otimes \cI) =L_A^2-c_2(F)= 9\) and, because 
    \[c_2(L \otimes \cI) = \sum_{p \in Z(\cI)} \mbox{dim}_{\bC} \cO_{A,p}/\cI_p,\]
    and \(\mbox{dim}_{\bC} \cO_{A,p}/\cI_p \geq 1\) for each \(p\), then \(Z(\cI)\) cannot have more than \(9\) points.
    Since the base locus of \(\varphi\) has at most \(9\) points and is invariant with respect to the action of \(K(L_B)\), it must be a translate of \(K(L_B)\). But, by \((-1)^{*}_BL_B \simeq L_B\) and the fact that there is an odd number of base points, the base locus must contain at least one \(2-\)torsion point.
\end{proof}
	
\begin{rem}
    Notice that the equality
    \[9 = c_2(L \otimes \cI) = \sum_{p \in Z(\cI)} \mbox{dim}_{\bC} \cO_{A,p}/\cI_p\]
    implies that \(\tilde{A}\) is the blow up of \(A\) once in each of the nine points of the base locus of \(\varphi\).
\end{rem}
	
We know that \(H^0(F)\) is generated by three sections \(s_1, s_2, s_3\) and that \(H^0(F)^{\vee} \cong \bigwedge^2 H^0(F)\) is generated by \(s_1 \wedge s_2, s_1 \wedge s_3, s_2 \wedge s_3\). Chosen a point \(p\) in the base locus of \(\varphi\), we can choose \(s_1, s_2, s_3\) such that
\[
\begin{cases}
    s_1(p) \neq 0,\\
    s_2(p) = s_3(p) = 0.
\end{cases}%
\]

Dualizing the exact sequence
\[
\begin{tikzcd}
0 \arrow{r} & F^\vee \arrow[r,"ev^*"] & H^0(F)^{\vee}\otimes \cO_A \arrow{r} & L_A\otimes \mathcal I \arrow{r} & 0 ,
\end{tikzcd}
\]
one sees that a point $p$ lies in the base locus $\mathrm{Spec}(\cO_A/\mathcal I)$ if and only if $F$ fails to be globally generated at $p$.
In particular, if $p$ is a base point then
\[
h^1\!\big(F\otimes \mathcal I_p\big)=1
\qquad\text{and}\qquad
h^0\!\big(F\otimes \mathcal I_p\big)=2,
\]
that is, imposing vanishing at $p$ does not cut out two independent linear conditions on $H^0(F)$. The curve \(Z(s_2 \wedge s_3) = \overline{ \bigcup_{\lambda \in \bC} Z(s_2 + \lambda s_3) } \) is hyperelliptic since the moving part of the vanishing locus of \(s_2 + \lambda s_3\) describes a \(g^1_2\) on \(H_p\). There exists such an hyperelliptic curve for any \(p\) in the base locus and they are one the translate of the other.\\
For any \(p\) in the base locus (with exceptional divisor \(E\)), we will denote \(H_p\) the associated hyperelliptic curve and \(H_E\) its strict transform. Let us remark that \(H_p\) is singular at \(p\) as it is the zero locus of the wedge of two sections vanishing at \(p\).

The following holds.

\begin{prop}\label{GenusHyperellCurves}
    The hyperelliptic curves associated to the exceptional divisors are either
    \begin{enumerate}
        \item smooth of genus \(6\) intersecting twice the associated exceptional divisor;
        \item nodal of genus \(5\) with one node in its only intersection with the associated exceptional divisor.
    \end{enumerate}
\end{prop}
	
\begin{proof}
    The map \(\varphi\) defined in \autoref{degree3cover} maps onto \(\bP(H^0(F)) = \bP(H^0(L_B)) \cong \bP^2\). The action of \(-1\) on \(H^0(L_B)\) has two non-trivial eigenspaces (see \cite[section 4.6]{BirkenhakeLange1992complex}), hence, the fixed point space of the map \(-1: \bP^2 \to \bP^2\) induced by \((-1)_A\) is the union of a line \(l\) and a point \(P\).\\
    Note that \((-1)_A\) extends to \(\tilde{A}\) and that the space \(\tilde{\varphi}^{*}(l) \cup \tilde{\varphi}^{*}(P)\) contains the fixed point space of \((-1)_{\tilde{A}}\).\\
    We note that the fixed point space of  \((-1)_{\tilde{A}}\) contains one and only one exceptional divisor \(E\). Since \(\tilde{\varphi}^{*}(P)\) consists in a finite number of points%
    , \(E \subset \tilde{\varphi}^{*}(l)\). But \(\tilde{\varphi}_{|_{\tilde{\varphi}^{*}(l)}}\) is a degree 3 map equivariant with respect to \(-1\), so the other irreducible component of \(\tilde{\varphi}^{*}(l)\) is a hyperelliptic curve \(\tilde{H}\), passing through at least $13$ $2$-torsion points (the fixed locus on $A$ minus $\tilde{\varphi}^{*}(P) $ ).
    \\
    This hyperelliptic curve is \(H_E\). Indeed, given a point in \(l \subset \bP H^0(F)\), it corresponds to the class \([s]\) of a section of \(H^0(F)\). Recall that \(\varphi^{-1}([s]) = Z(s)\). Since one of the preimages of \([s]\) via \(\tilde{\varphi}\) is in \(E\), one of its zeroes is the corresponding base point of \(\varphi\). This means that every \([s] \in l\) vanishes in the base point \(p\) corresponding to the exceptional divisor \(E\), that is, every element of \(l\) is the class of a linear combination of the sections \(s_1, s_2\) such that \(Z(s_1 \wedge s_2) = H_p\). Since \(\tilde{\varphi}(\tilde{H}) = \varphi(H_p)\), \(\tilde{H}\) is the strict transform of \(H_p\) (that is \(H_E\)).\\    
    Now, since \(H_p\) lies in \(|L_A|\), \(p_a(H_p) = 7\). Therefore, \(p_g(H_p) \leq 7\). Moreover, \(\varphi_{|_{H_p}}\) is the degree \(2\) map \(H_p \to \bP^1 \cong l\) and its ramification contains \(A[2] \smallsetminus \tilde{\varphi}^{*}(P)\). %
   
    Let $w=2p_g(H_p)+2$ denote the number of Weierstrass points of $H_p$. Call $s,n,c,d$ respectively the number of smooth points, nodes, cusps and  singularities with $\delta $ invariant $\ge 2$,  among $H_p\cap A[2]$, in this way we have  $13\le s+n+c+d$. Then, 
each of the $2-$torsion points in $H_p\cap A[2]$ contributes the number of Weierstrass points of (a desingularization of) $H_p$ in the following way:
\begin{itemize}
    \item If $p\in (H_p\cap A[2])$ is a smooth point or a cusp for $H_p$, then it induces one Weierstrass point on the desingularization of $H_p$;
\item If $p\in (H_p\cap A[2])$ is a node for $H_p$, then it induces two Weierstrass points on the desingularization of $H_p$.
\end{itemize}
 Therefore, we get $
    w\ge 2n+c+s\ge   13+n-d$ (with equality holding if $d=0$ and $|H_p\cap A[2]|=13$) and $p_g\le 7-n-c-2d$ (in this case the equality implies that $H_p$ is smooth away from the $2-$torsion points and that $\delta=-2$ for all points with a singularity worse than a node or a cusp). Putting all together, we get 
    \[
    14-2n-2c-4d+2\ge 13+n-d.
    \]
    Hence $3-3n-2c-3d\ge 0.$ Taking into account that $n,c,d\ge 0$ and $n+d+c\ge 1$ we obtain $2$ possibilities:
    \begin{enumerate}
        \item \(H_E\) is a smooth hyperelliptic curve of genus \(6\) corresponding to $n=1,c=0,d=0$;
        \item \(H_E\) is a hyperelliptic curve of genus \(5\) with one node in its intersection with the exceptional divisor, corresponding to $n=0,c=0,d=1$ (In this case, notice that we must have that the delta invariant is exactly $-2$ since $w=12$ and the only such singularity for $H_p$ making the number of Weiestrass points $w=13-d=12$ is $y^2=x^4)$.
    \end{enumerate}
    We cannot have $c=1,d=0,n=0$ because otherwise we would have an odd $w$ which is impossible.
\end{proof}

\begin{figure}[ht]
    \centering
    \includegraphics[width=0.7\linewidth]{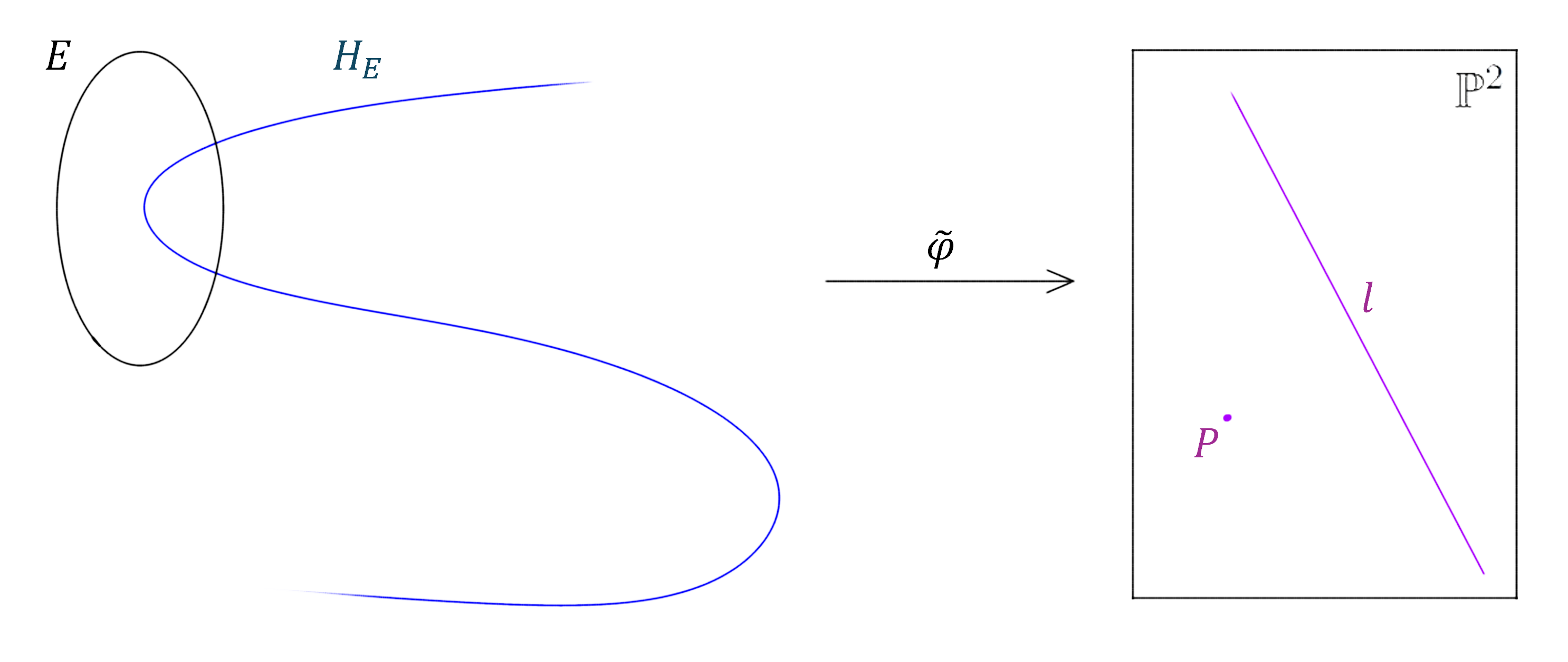}
    \caption{\(\tilde{\varphi}\) maps \(E\) and \(H_E\) to the line \(l\) in the fixed point space of \(-1 \colon \bP^2 \to \bP^2\).}
    \label{fig:FixPtsLocusOpposite}
\end{figure}

\section{The case of \texorpdfstring{$B = \EllCurve \times \EllCurve$}{product of elliptic curves}}

\label{Section: 3}
We would like to analyze more in detail the case of \((1,3)\) abelian surfaces of the form \(B = \EllCurve\times \EllCurve\) (with \(\EllCurve\) an elliptic curve) endowed with the polarization $L_B=\cO_{\EllCurve \times \EllCurve}(\EllCurve_1+\EllCurve_2+\Delta)$ where $\EllCurve_1=\EllCurve \times \{0\}$, $\EllCurve_2= \{0\} \times \EllCurve$ and $\Delta$ denotes the diagonal. Let $b_1\neq b_2\in \EllCurve[2]$ be non trivial distinct $2$ torsion points and $A=\EllCurve\times \EllCurve/\langle (b_1,b_2) \rangle$ be the associated $(1,6)$ abelian surface.\\
To the purpose of this work, this is more than a simple example, indeed, we will prove that, in these conditions, the hyperelliptic curves associated to the exceptional divisors are smooth of genus \(6\); therefore, this is the case for the general \((1,6)\) abelian surface.\\
First let us compute the kernel of the polarization $L_B$.

\begin{lem}
    The kernel of the polarization \(L_B\) is \(K(L_B)=\{(x,-x) | x\in \EllCurve[3]\}\).
    \label{kernel}
\end{lem}

\begin{proof}
    We consider the family of decomposable divisors linearly equivalent to \(\EllCurve_1 + \EllCurve_2 + \Delta\) obtained translating \(\EllCurve_1\), \(\EllCurve_2\) and \(\Delta\) independently; when the translations are given by the same point for all three components, said point is in the kernel.\\
    First we notice that, for all \((x,y) \in \EllCurve \times \EllCurve\), 
    \[t_{(x,y)}\EllCurve_1 = t_{(0,y)}\EllCurve_1, \quad t_{(x,y)}\EllCurve_2 = t_{(x,0)}\EllCurve_2, \quad t_{(x,y)}\Delta = t_{(\frac{x-y}{2},-\frac{x-y}{2})}\Delta,\]
    that is each of these divisors is not affected by translations in its own direction. So, we want the conditions in which
    \[t_{(0,y)}\EllCurve_1 + t_{(x,0)}\EllCurve_2 + t_{(a, -a)}\Delta \sim \EllCurve_1 + \EllCurve_2 + \Delta.\]
    By the theorem of the cube \(t_{(0,y)}\EllCurve_1 + t_{(x,0)}\EllCurve_2 + t_{(a, -a)}\Delta - (\EllCurve_1 + \EllCurve_2 + \Delta) \sim 0\) if and only if its restrictions to all the translations of \(\EllCurve_1\) and \(\EllCurve_2\) are linearly equivalent to \(0\). This happens when \(x = -2a\) and \(y = 2a\); so, for all \(a \in \EllCurve\), \(t_{(0,2a)}\EllCurve_1 + t_{(-2a,0)}\EllCurve_2 + t_{(a, -a)}\Delta\) is linearly equivalent to \(\EllCurve_1 + \EllCurve_2 + \Delta\), and these are all the divisors of this kind.\\
    Finally, a point \((x,y) \in \EllCurve \times \EllCurve\) is in the kernel \(K(L_B)\) if and only if there exist \(a \in \EllCurve\) such that
    \[
    \begin{cases}
        x = -2a = a\\
        y = 2a = -a
    \end{cases}%
    \]
    that is, \(K(L_B)=\{(x,-x) | x\in \EllCurve[3]\}\).
\end{proof}

A simple observation is that the chain of maps $\EllCurve\times \EllCurve\to A \dashrightarrow \mathbb P^2$ is equivariant with respect to $K(L_B)$, hence the base locus of $\EllCurve\times \EllCurve \dashrightarrow \mathbb P^2$ are $18$ points invariant under translations by $K(L_B)$. We deduce that an appropriate translation of the antidiagonal $\Delta^-_a=\{(a+x,a-x)\quad x\in \EllCurve\}\subset \EllCurve\times \EllCurve$ contains $9$ points of the base locus (a translate of $K(L_B)$). Any translate of $\Delta^-$ is equivariant with respect to the translations by $K(L_B)$. We get that the induced map $\EllCurve \simeq \Delta^-_a\to \mathbb P^2$ is of degree $3$ and equivariant with respect to $K(L_B)\simeq \EllCurve[3]$. We deduce readily.

\begin{lem}
    The action of $K(L_B)$ on $\mathbb P^2$ is (after fixing coordinates) generated by the matrices
    \[
    M=\begin{pmatrix}
        1 & 0 & 0 \\
        0 & \zeta_3 & 0 \\
        0 & 0 & \zeta_3^2
    \end{pmatrix},%
    \quad%
    N=\begin{pmatrix}
        0 & 0 & 1 \\
        1 & 0 & 0 \\
        0 & 1 & 0
    \end{pmatrix},
    \] 
    where $\zeta_3^3=1.$
\end{lem}

\begin{proof}
    After a suitable choice of coordinates, we can assume that \(\tilde{\varphi}\) maps the strict transform of the antidiagonal \(\Delta^-\) to the curve
		\[C_a = \{(x : y : z) \in \bP^2 | \, x^3 + y^3 + z^3+axyz = 0\},\]
	where \(a \in \mathbb{C}\) depends on the isomorphism class of \(C_a\).
	
	Since the base points of the pencil are the flexes of \(C_a\), the action of \(K(L_B)\) has to permute them. This condition is satisfied only by the group generated by the matrices \(M\) and \(N\) in the statement. Therefore, \(M\) and \(N\) generate the representation of \(K(L_B)\) on \(\bP^2\).
\end{proof}

Employing the family of decomposable divisors defined in \autoref{kernel}, we can find the base locus of \(\varphi\).

\begin{lem}
    The base locus of \(\varphi\) is \([(b_1 + b_2, b_1)] + K(L_B)\).
    \label{baselocus}
\end{lem}

\begin{proof}
    Note that, given \(\sigma \in H^0(\pi_* \cO(\EllCurve_1 + \EllCurve_2 + \Delta))\), there exists \(s \in H^0(\cO(\EllCurve_1 + \EllCurve_2 + \Delta))\) such that, around \([p] \in A = \EllCurve \times \EllCurve / \langle b \rangle\) (with \(b \coloneqq (b_1, b_2)\)),
    \[\sigma([p]) =
    \begin{pmatrix}
        s(p)\\
        s(p + b)
    \end{pmatrix}.
    \]
    Then, \([p] \in Z(\sigma)\) if and only if \(p \in Z(s) \cap Z(t_b^* s)\). When \(Z(s) = t_{(0,2a)}\EllCurve_1 + t_{(-2a,0)}\EllCurve_2 + t_{(a, -a)}\Delta\), then
    \[
    Z(s) \cap Z(t_b^* s) = \left\{ 
    \begin{aligned}
        &(-2a + b_1, 2a),     &   &(4a + b_1 + b_2, 2a),   &   &(-2a, 4a + b_1 + b_2),\\ 
        &(-2a, 2a + b_2),   &   &(4a + b_2, 2a + b_2),   &   &(-2a + b_1, -4a + b_1)
    \end{aligned} \right\}
    \]
    and we find that:
    \begin{align*}
        &(b_1 + b_2, b_1) = (-2a, 2a + b_2)    &   &\iff   &   &2a = b_1 + b_2\\
        &(b_1 + b_2, b_1) = (4a + b_1 + b_2, 2a)   &   &\iff    &     &2a = b_1\\
        &(b_1 + b_2, b_1) = (-2a + b_1, -4a + b_1)     &    &\iff  &   &2a = b_2
    \end{align*}
    this means that there are three sections (linearly independent in pairs) in the linear system of \(\EllCurve_1 + \EllCurve_2 + \Delta\) that vanish in \((b_1 + b_2, b_1)\). This implies the statement.  
\end{proof}

\begin{prop}
    The hyperelliptic curves associated with the exceptional divisors are smooth of genus \(6\).
    \label{node}
    \end{prop}

\begin{proof}
    Called \(\sigma\), \(\rho\) and \(\tau\) the sections found in \autoref{baselocus} respectively correspondent to \(2a = b_1 + b_2, b_1, b_2\), and fixed coordinates \(x\) and \(y\) around \((b_1 + b_2, b_1)\) that vanish respectively on the translates of \(\EllCurve_2\) and \(\EllCurve_1\), we find the following local equations:
    \begin{align*}
        \sigma(x,y) = x + O(x^2,xy,y^2), \quad &t^*_b \sigma(x,y) = \mu y + O(x^2,xy,y^2),\\
        \rho(x,y) = \lambda y + O(x^2,xy,y^2), \quad &t^*_b \rho(x,y) = x - y + O(x^2,xy,y^2),\\
        \tau(x,y) = x - y + O(x^2,xy,y^2), \quad &t^*_b \tau(x,y) = \nu x + O(x^2,xy,y^2)
    \end{align*}
    where \(\lambda\), \(\mu\) and \(\nu\) are complex numbers; so, around \([(b_1 + b_2, b_1)]\),
    \begin{align*}
        \pi_{*}\sigma(x,y) &= 
        \begin{pmatrix}
            x + O(x^2,xy,y^2)\\
            \mu y + O(x^2,xy,y^2)
        \end{pmatrix},\\
        \pi_{*}\rho(x,y) &= 
        \begin{pmatrix}
            \lambda y + O(x^2,xy,y^2)\\
            x - y + O(x^2,xy,y^2)
        \end{pmatrix},\\
        \pi_{*}\tau(x,y) &=
        \begin{pmatrix}
            x - y + O(x^2,xy,y^2)\\
            \nu x + O(x^2,xy,y^2)
        \end{pmatrix}.
    \end{align*}
    Since they are not linearly independent there exist \(\alpha, \beta \in \bC\) such that \(\pi_{*}\sigma = \alpha \pi_{*}\rho + \beta \pi_{*}\tau\); confronting the coefficients of the coordinates, one finds \(\lambda \nu = -1\) and, so, 
    \[\pi_{*}\rho \wedge \pi_{*}\tau = xy - x^2 - y^2,\]
    which is not a square. Thus, the hyperelliptic curve associated with \([(b_1 + b_2, b_1)]\) has a node in \([(b_1 + b_2, b_1)]\), and its strict transform is a smooth hyperelliptic curve of genus \(6\).
\end{proof}
As a corollary, we obtain \autoref{introhyp}. 

\begin{cor}
    On a very general $(1,6)$ abelian surface $(S,L)$ there is a unique hyperelliptic curve of geometric genus $6$ in $|L|$.
\end{cor}

\begin{proof}
    The existence part follows from \autoref{node}. For uniqueness, consider a hyperelliptic curve of geometric genus $6$ singular at $p$ (the curve is singular at a unique point, otherwise the geometric genus would be $<6$), together with $A$ the $g^1_2$ (defined on the desingularization $\widetilde{D})$. We denote with $i:\widetilde{D}\to A$ the natural map. We get an associated kernel bundle
    \[
    \begin{tikzcd}
        0 \arrow{r} & E_{D,A}^\vee \arrow{r} & H^0(A) \otimes \mathcal{O}_S \arrow{r}{ev} & i_*A.
    \end{tikzcd}
    \]
    If $\mathrm{NS}(S)=\mathbb Z\cdot L$, then $E_{D,A}$ is stable with $c_2=2,3$. Indeed, the sections $s\in H^0(A)^\vee\subset H^0(E_{D,A})$ vanish at the preimages of the $g^1_2$ and with order at most one at $p$ (if the sections were vanishing at higher order at $p$ then $D=Z(\bigwedge ^2 H^0(A)^\vee)$ would have a worst singularity then a node or a cusp at $p$, hence it would be of geometric genus $<6$). Now $c_2$ cannot be $2$ by stability reasons, hence $c_2=3$. This means that $E_{D,A}$ is the semihomogeneous vector bundle with primitive first Chern class, rank $2$ and $c_2=3$ (which is unique up to twists by topologically trivial line bundles, c.f. \cite{MUKAI}). Hence $H^0(A)^\vee=H^0(E \otimes     \mathcal{I}_p)$ and $D=H_p$ is the hyperelliptic curve associated with $p$ (i.e. $D=H_p=Z(\bigwedge ^2 H^0(E \otimes     \mathcal{I}_p))$ and $\widetilde{D}=H_E$).
    \end{proof}

\section{Ramification of the degree \texorpdfstring{\(3\)}{3} cover}
\label{Section: 4}
\label{sec:ramification}
In the following section we will investigate the nature of the ramification of \(\tilde{\varphi} \colon \tilde{A} \to \bP^2\), the resolution of singularities of the rational map \(\varphi\) defined in \autoref{degree3cover}. In order to do so, we notice that the projection \(\pi \colon B \to A\) and the morphism \(\varphi\) are equivariant with respect of the action of \(K(L_B)\) which, to simplify notation, we will call \(\ShortKer\) from now on. Hence, taking the quotient by \(\ShortKer  \cong \bZ_3 \times \bZ_3\) and by any of its proper subgroups \(\Subgr \cong \bZ_3\), we get the diagram
\begin{equation}\label{eq:DiagramQuotientsKLB}
\begin{tikzcd}
    B \arrow[r, "\pi"] \arrow[d] & A \arrow[r, dashed, "\varphi"] \arrow[d, "p_1"] & \bP^2 \arrow[d]\\
    B/\Subgr \arrow[r] \arrow[d] & A/\Subgr \arrow[r, dashed] \arrow[d, "p_2"] & \bP^2/\Subgr \arrow[d]\\
    B/\ShortKer \arrow[r] & A/\ShortKer \arrow[r, dashed,"\varphi'"] & \bP^2/\ShortKer 
\end{tikzcd}
\end{equation}
where the vertical arrows of the first two columns are \'{e}tale covers.

We will show that \(A/\Subgr\) is a \((1,2)\) abelian surface, and \(A/\ShortKer\) is a \((1,6)\) abelian surface.\\
Notice that the action of \(\ShortKer\) over \(A\) extends naturally to \(\tilde{A}\), mapping the exceptional divisors to each other. This means that \(\tilde{A}/\ShortKer\) is the blow up in one point of a \((1,6)\) abelian surface and the morphism \(\tilde{\varphi}' \colon \tilde{A}/\ShortKer \to \bP^2/\ShortKer \) is the resolution of singularities of \(\varphi' \colon A/\ShortKer \dashrightarrow \bP^2/\ShortKer\).\\
Throughout this section, we will call:
\begin{itemize}
    \item \(p \coloneqq p_2 \circ p_1\);
    \item \(\tilde{p} \colon \tilde{A} \to \tilde{A}/\ShortKer\) the projection;
    \item \(\epsilon \colon \tilde{A} \to A\) the blow up of the surface at the nine points of the base locus.
\end{itemize}
The purpose of this section is to prove that the image of the ramification locus of \(\tilde{\varphi}\) via \(\tilde{p}\) (that is, the quotient of the ramification by the action of \(\ShortKer\)) is a genus 6 smooth hyperelliptic curve.

\begin{figure}[ht]
    \centering
    \includegraphics[width=0.8\linewidth]{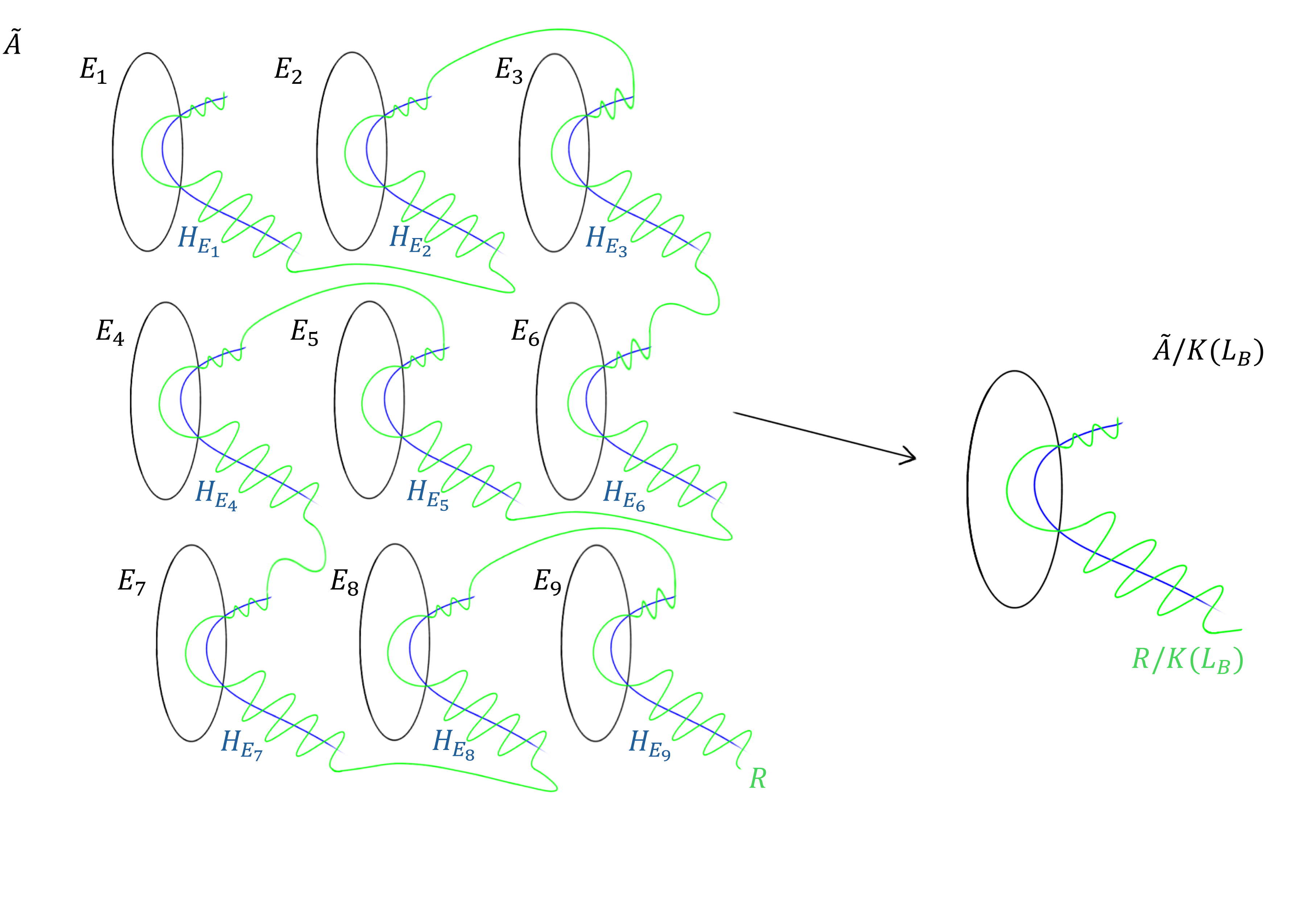}
    \caption{The quotient by \(\ShortKer = K(L_B)\) on the exceptional locus, the hyperelliptic curves and the ramification of \(\tilde{\varphi}\).%
    }
    \label{fig:ProiezioneGrado9}
\end{figure}

\begin{prop}
    Let \(R\) be the ramification of \(\tilde{\varphi}\) and let \(E_1, \dots, E_9\) be the exceptional divisors of \(\tilde{A}\). Then, \(\cO(R) = \cO_{\tilde{A}}(3) \otimes \cO_{\tilde{A}}\left(\sum_{i=1}^9 E_i\right)\). In particular,
    \[p_a(R) = 46, \quad R^2 = 72\]
\end{prop}

\begin{proof}
    Since \(\tilde{\varphi}\) is a finite morphism between projective varieties, we get by Hurwitz formula
    \[K_{\tilde{A}} \sim \tilde{\varphi}^{*} K_{\bP^2} + R,\]
    where \(\sim\) denotes the linear equivalence.  Moreover, taking into account \(\cO_{\tilde{A}}(K_{\tilde{A}}) = \cO_{\tilde{A}}(\sum_{i=1}^9 E_i)\) and \(\cO_{\bP^2}(K_{\bP^2}) = \cO(-3)\), we get:
    \[\cO_{\tilde{A}}(R) = \cO_{\tilde{A}}(3) \otimes \cO_{\tilde{A}}\left(\sum_{i=1}^9 E_i\right).\]
    It follows \(R^2 = 72\). Furthermore, by Riemann-Hurwitz,
    \[2p_a(R) - 2 = R(R + K_{\tilde{A}}) = 90,\]
    hence, \(p_a(R) = 46\).
\end{proof}

\begin{prop}
    With the notation above, called \(\Subgr \cong \bZ_3\) a proper nontrivial subgroup of \(\ShortKer \cong \bZ_3 \times \bZ_3\):
    \begin{enumerate}
        \item \(A/\Subgr\) is a \((1,2)\) abelian surface,
        \item \(A/\ShortKer\) is a \((1,6)\) abelian surface.
    \end{enumerate}
\end{prop}

\begin{proof}
\noindent
    \begin{enumerate}
        \item If \(A = \bC^2 / \Lambda_A\), there exist a basis \(\{\lambda_1, \lambda_2, \mu_1, \mu_2\}\) of \(\Lambda_A\) such that the  alternating 2-form \(E\) associated to the polarization \(L_A\) is represented by the matrix
        \[\begin{pmatrix}
            0 & 0 & 1 & 0\\
            0 & 0 & 0 & 6\\
            -1 & 0 & 0 & 0\\
            0 & -6 & 0 & 0\\
        \end{pmatrix}.\]
        \\Moreover, \(K(L_A) = \ker(E) / \Lambda_A\), where
        \[\ker(E) = \{v \in \bC^2 | E(v, \lambda) \in \bZ \; \forall \lambda \in \Lambda_A\} = \Bigl \langle \lambda_1, \frac{1}{6} \lambda_2, \mu_1, \frac{1}{6} \mu_2 \Bigr \rangle_{\bZ},\]
        hence \(\ShortKer = K(L_B) = \Bigl \langle \lambda_1, \frac{1}{3} \lambda_2, \mu_1, \frac{1}{3} \mu_2 \Bigr \rangle_{\bZ} / \Lambda_A\), and \(A/G\) is \(\bC^2 / \Lambda'\) where \(\Lambda'\) is either \(\Bigl \langle \lambda_1, \frac{1}{3} \lambda_2, \mu_1, \mu_2 \Bigr \rangle_{\bZ}\), \(\Bigl \langle \lambda_1, \lambda_2, \mu_1, \frac{1}{3} \mu_2 \Bigr \rangle_{\bZ}\), \(\Bigl \langle \lambda_1, \frac{1}{3} \lambda_2 + \frac{1}{3} \mu_2, \mu_1, \mu_2 \Bigr \rangle_{\bZ}\) or \mbox{\(\Bigl \langle \lambda_1, \lambda_2, \mu_1, \frac{1}{3} \lambda_2 + \frac{1}{3} \mu_2 \Bigr \rangle_{\bZ}\)} depending on which proper nontrivial subgroup of \(\ShortKer = K(L_B)\) is considered. With respect to these bases, \(E\) is represented by the matrix
        \[\begin{pmatrix}
            0 & 0 & \mp 1 & 0\\
            0 & 0 & 0 & \mp 2\\
            \pm 1 & 0 & 0 & 0\\
            0 & \pm 2 & 0 & 0\\
        \end{pmatrix},\]
        that is, \(A/G\) is a \((1, 2)\) abelian surface.
        \item Since the ramification lies in \(|3 \cO(1) + K_{\tilde{A}}| = |3 \epsilon^{*}L_A -2 \sum_{i=1}^9 E_i|\), \(\epsilon_{*}R\) lies in \(|3 L_A|\), so \((\epsilon_{*}R)^2 = 9 \cdot L_A^2 = 9 \cdot 12\). Then, \((p_{*}\epsilon_{*}R)^2 = \frac{(\epsilon_{*}R)^2}{9} = 12\), hence it is a \((1,6)\) polarization for \(A/\ShortKer\).
    \end{enumerate}
\end{proof}

\begin{prop}\label{QuotientRamification}
    \(\tilde{p}(R)\) is a genus 6 smooth hyperelliptic curve.
\end{prop}

\begin{proof}
    Since \(\varphi\) is equivariant with respect to \((-1)_A\), we notice that the ramification \(R\) is invariant with respect to \((-1)_{\tilde{A}}\). Moreover, since, as seen in \autoref{GenusHyperellCurves}, the degree \(2\) map on the hyperelliptic curve \(H_E\) associated to the exceptional divisor \(E\) is \(\tilde{\varphi}_{|_{H_E}}\), the ramification points of \(H_E\) are ramification points for \(\tilde{\varphi}\). Thus, \(\epsilon(R)\) contains \(13\) \(2\)-torsion points (the ramification points of \(H_E\)); then \(\epsilon(R)/\ShortKer\) contains \(13\) \(2\)-torsion points of \(A/\ShortKer\) and its strict transform \(\tilde{p}(R) = R / \ShortKer\), that is the ramification of the morphism \(\tilde{A}/\ShortKer \to \bP^2/\ShortKer\), has \(14\) fixed points with respect to \((-1)_{\tilde{A}/\ShortKer}\).\\
    Now, the quotient $\tilde{p}(R) \to \tilde{p}(R)/(-1)$ is a degree 2 rational map.   As in the proof of \autoref{GenusHyperellCurves}, let   \(n\) be the number of nodes, \(c\) the number of cusps and \(d\) the number of singularities with \(\delta\) invariant \(\geq 2\) of \(\tilde{p}(R)\). The resolution of this degree 2 map has a number of ramification points greater or equal to \(14 + n - d\). Moreover, applying the Riemann-Hurwitz formula to the degree \(9\) \'{e}tale cover \(R \to \tilde{p}(R)\), we get \(g(\tilde{p}(R)) \leq 6 - n - c - 2d\).\\
    Combining these inequalities with the Riemann-Hurwitz formula for \(\tilde{p}(R) \to \tilde{p}(R)/(-1)_{\tilde{A}/\ShortKer}\), we get
    \[2(6 - n - c - 2d) - 2 \geq 2 (2g(\tilde{p}(R) / (-1)) - 2) + 14 + n - d.\]
    Hence \(-3n -2c - 3d - 4 g(\tilde{p}(R) / (-1)) \geq 0\). Taking into account that \(n, c, d, g(\tilde{p}(R)/ (-1)) \ge 0\), the only possible solution is \(n = c = d = g(\tilde{p}(R)/ (-1)) = 0\), that is, \(\tilde{p}(R)/ (-1)\) is the projective line and \(\tilde{p}(R)\)  is a genus 6 smooth hyperelliptic curve.
\end{proof}

\begin{cor}
    The ramification of \(\tilde{\varphi}\) is smooth of geometric genus \(46\).
\end{cor}

\section{Invariants of the Galois closure of \texorpdfstring{\(\varphi\)}{the degree 3 cover}}
\label{Section: 5}
We recall briefly a geometrical description of the Galois closure \(\GalCl\) of a degree \(\degrm\) generically finite morphism \(\RatMap \colon Z \to Y\) as reported on \cite{BPS}.

Chosen \(\cW \subset Y\) a non-empty Zariski open subset of \(Y\) such that \(\RatMap_{|_{\RatMap^{-1}(\cW)}} \colon \RatMap^{-1}(\cW) \to \cW\) is an \'{e}tale morphism, denote
\[\GalCl^{\circ} \coloneqq \left\{(z_1, \dots, z_{\degrm}) \in Z^\degrm \;|\; \exists y \in \cW \mbox{ such that } \{z_1, \dots, z_{\degrm}\} = \RatMap^{-1}(y)\right\}.\]
By lifting arcs, the irreducible components of \(\GalCl^{\circ}\) are all isomorphic. The normalization \(\GalCl\) of the Zariski closure of any fixed irreducible component of \(\GalCl^{\circ}\) together with \(\delta\), the composition of the \(i\)-th projection \(\pi_i \colon \GalCl \to Z\) and \(\RatMap\), is the Galois closure of \(\RatMap\).
\[
\begin{tikzcd}
    W \arrow[d, "\pi_i"] \arrow[dr]\\
    Z \arrow[r, "\gamma"] & Y
\end{tikzcd}
\]

\begin{rem}
    We note that the variety \(\GalCl^{\circ}\) can also be defined as a subvariety of \(Z^{\degrm - 1}\)
    \[\GalCl^{\circ} \coloneqq \left\{(z_1, \dots, z_{\degrm - 1}) \in Z^\degrm \;|\; \exists y \in \cW \mbox{ such that } \exists z_{\degrm} \in \RatMap^{-1}(\cW) \mbox{ such that } \{z_1, \dots, z_{\degrm}\} = \RatMap^{-1}(y)\right\}.\]
\end{rem}

We study the Galois closure of the morphism \(\tilde{\varphi}\). Denote \(\GalCl\) the Zariski closure of 
\[\tilde{A} \times_{\bP^2} \tilde{A} -\Delta= \left\{(p,q) \in \tilde{A} \times \tilde{A} \;|\; \tilde{\varphi}(p) = \tilde{\varphi}(q) ,p\neq q\right\}.\]
We will show that it is smooth, and, hence, it is the Galois closure of \(\tilde{\varphi}\).

The first projection \(\tilde{A} \times_{\bP^2} \tilde{A} \to \tilde{A}\) induces a degree \(2\) map \(\delta \colon \GalCl \to \tilde{A}\) such that, given \(s \in \bP^2\), the cardinality of its preimage determines if the preimages are in \(B_{\delta}\), the branch locus of \(\delta\).
\begin{itemize}
    \item If \(|\tilde{\varphi}^{-1}(s)| = 3\), then \(\tilde{\varphi}^{-1}(s)\) consists in three distinct points \(\{p_1, p_2, p_3\}\) and, for \(i = 1, 2 ,3\), one gets \(\delta^{-1}(p_i) = \{(p_i, p_j)\}_{j \in \{1, 2, 3\} \smallsetminus \{i\}}\), so \(p_i\) is not in the branch locus of \(\delta\);
    \item if \(|\tilde{\varphi}^{-1}(s)| = 2\), then \(\tilde{\varphi}^{-1}(s)\) consists in one point \(p\) with multiplicity \(1\) for \(\tilde{\varphi}\) and one point \(q\) with multiplicity \(2\) for \(\tilde{\varphi}\); then \(\delta^{-1}(q) = \{(q, p), (q, q)\}\) and \(\delta^{-1}(p) = \{(p,q)\}\), that is, \(p\) is a branch point for \(\delta\);
    \item if \(|\tilde{\varphi}^{-1}(s)| = 3\), the only point in \(\tilde{\varphi}^{-1}(s)\) is a triple point for \(\tilde{\varphi}\) and, hence, a branch point for \(\delta\).
\end{itemize}
Therefore, set-wise, \(B_{\delta}\) consists of the triple points of \(\tilde{\varphi}\) and of the single points points of \(\tilde{\varphi}\) \(p \in \tilde{A}\) such that  \(\tilde{\varphi}^{-1}(\tilde{\varphi}(p))\) contains a double point for \(\tilde{\varphi}\).

\begin{prop}
    The branch locus \(B_{\delta}\) of \(\delta\) is a smooth curve with
    \[p_a(B_{\delta}) = 46, \quad B_{\delta}^2 = -36\]
\end{prop}

\begin{proof}
    If we call \(B\) the branch divisor of \(\tilde{\varphi}\) and \(R\) the ramification divisor of \(\tilde{\varphi}\), by what we noticed above, we get
    \[\tilde{\varphi}^{*}(B) = B_{\delta} + 2R\]
    and, since \(R \sim 3 \epsilon^{*}L_A -2 \sum_{i=1}^9 E_i\) and \(B \sim \tilde{\varphi}_{*}R\), that is \(\cO(B) = \tilde{\varphi}_{*} (\cO(\sum_{i=1}^9 E_i) \otimes \tilde{\varphi}^{*}\cO(3)) = \cO(9) \otimes \deg\tilde{\varphi} \cO(3) = \cO(18)\), we get \(\tilde{\varphi}^{*}B \sim 18 \epsilon^{*}L_A - 18 \sum_{i=1}^9 E_i\) and
    \[B_{\delta} \sim 12\epsilon^{*}L_A - 14 \sum_{i=1}^9 E_i,\]
    from which, we readily get
    \[2 p_a(B_{\delta}) - 2 = (12\epsilon^{*}L_A - 14 \sum_{i=1}^9 E_i)(12\epsilon^{*}L_A - 13 \sum_{i=1}^9 E_i) = 90,\]
    that is \(p_a(B_{\delta}) = 46\), and
    \[B_{\delta}^2 = (12\epsilon^{*}L_A - 14 \sum_{i=1}^9 E_i)^2 = -36.\]
    Now, we notice that the map \(\eta \colon R \to B_{\delta}\) mapping a double point \(q\) of \(\tilde{\varphi}\) to the other point in the preimage of \(\tilde{\varphi}(q)\) and fixing the triple points of \(\tilde{\varphi}\) is birational. Hence \(p_g(B_{\delta}) = p_g(R) = 46\), and \(B_{\delta}\) is smooth (since the arithmetic genus and the geometric genus coincide). 
\end{proof}

\begin{figure}[ht]
    \centering
    \includegraphics[width=0.8\linewidth]{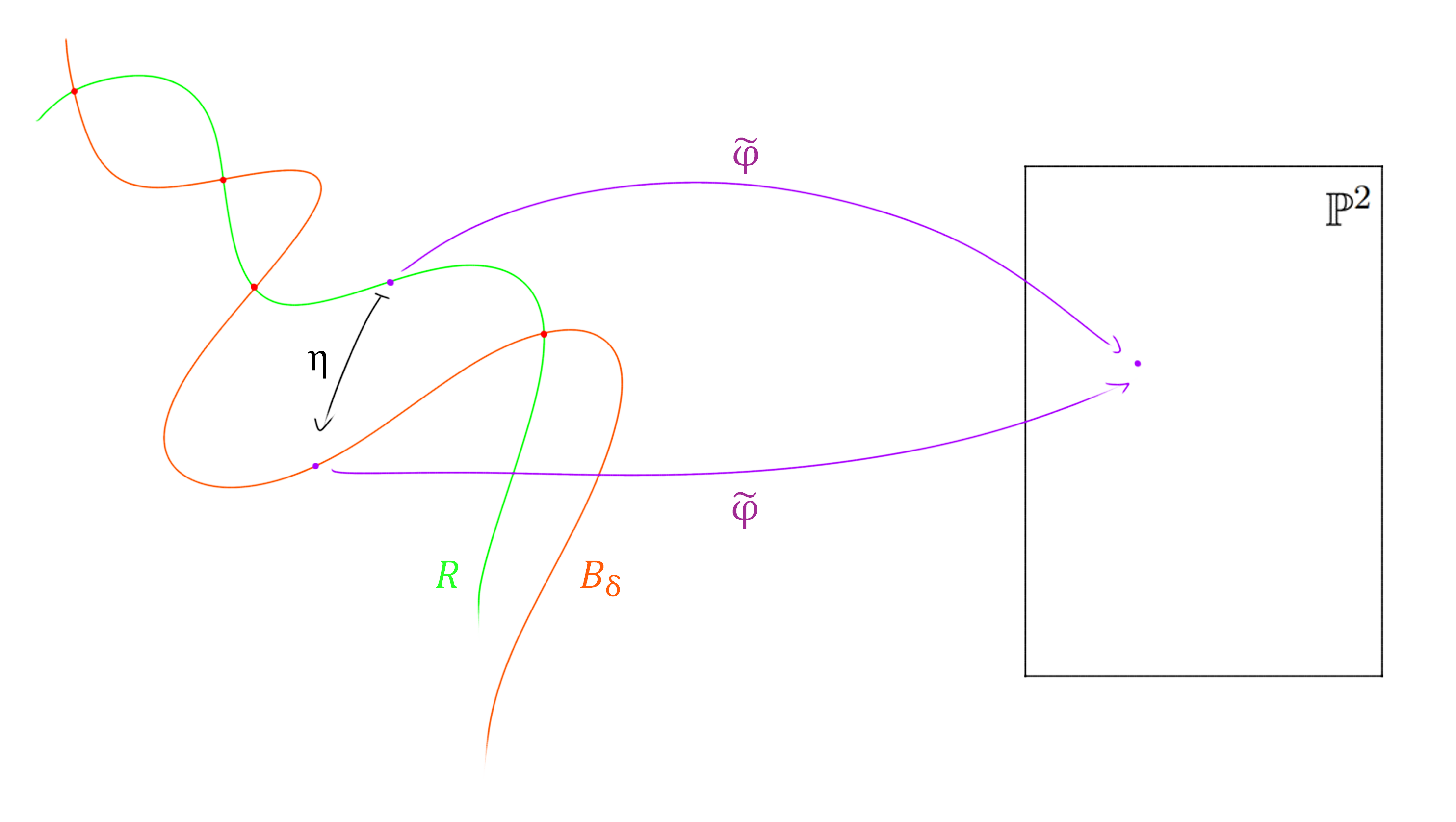}
    \caption{The map between the ramification of \(\tilde{\varphi}\) and the branch locus of \(\delta\).}
    \label{fig:MappaBranchRamification}
\end{figure}

\begin{rem}\label{rem:Wsmooth}
    The Galois closure of a degree $3$ finite morphism is smooth as long as the ramification and the target are smooth. Hence the surface \(W\) is embedded in \(\tilde{A} \times \tilde{A}\) because it is the Galois closure of \(\tilde{\varphi}\). 
\end{rem}

\begin{cor}
    \(\GalCl\) is smooth and, hence, it is the Galois closure of \(\tilde{\varphi}\).
\end{cor}


\begin{thm}
    The Galois closure \(\GalCl\) of \(\tilde{\varphi}\) has invariants:
    \[K_{\GalCl}^2 = 216, \quad c_2(\GalCl) = 108, \quad \chi(\cO_{\GalCl}) = 27.\]
\end{thm}

\begin{proof}
    Similarly to the case studied in \cite{BPS}, following \cite{persson}, the formulas to compute the invariants of \(\GalCl\) are:
    \[K_{\GalCl}^2 = 2(K_{\tilde{A}}^2 + 2 p_a(B_{\delta}) - 2) - \frac{3}{2} B_{\delta}^2,\] 
    \[c_2(\GalCl) = 2 c_2(\tilde{A}) + 2 p_a(B_{\delta}) - 2,\] 
    \[\chi(\cO_{\GalCl}) = 2 \chi(\cO_{\tilde{A}}) + \frac{p_a(B_{\delta}) - 1}{2} - \frac{B_{\delta}^2}{8}.\]
    Noticing that \(\tilde{A}\) is obtained by blowing up \(9\) times the abelian surface \(A\), and so \(c_2(\tilde{A}) = 9\) and \(\chi(\cO_{\tilde{A}}) = 0\), we get the values reported in the statement.
\end{proof}

We, now, notice that the action of \(\ShortKer = K(L_B)\) extends naturally to \(\GalCl\), the Galois closure of \(\tilde{\varphi}\). The action of \(\ShortKer\) on \(\GalCl\) is the restriction of the action of \(\{(k,k) \in \tilde{A} \times \tilde{A} | k \in \ShortKer\} \cong \ShortKer\) on \(\tilde{A} \times \tilde{A}\). Let \(\Subgr\) be a proper subgroup of \(\ShortKer\). Since the quotient maps \(\GalCl \to \GalCl/\Subgr\) and \(\GalCl \to \GalCl/\ShortKer\) are \'{e}tale coverings, \(\GalCl/\Subgr\) and \(\GalCl/\ShortKer\) have invariants:
\[K_{\GalCl/\Subgr}^2 = 72, \quad c_2(\GalCl/\Subgr) = 36, \quad \chi(\cO_{\GalCl/\Subgr}) = 9;\]
\[K_{\GalCl/\ShortKer}^2 = 24, \quad c_2(\GalCl/\ShortKer) = 12, \quad \chi(\cO_{\GalCl/\ShortKer}) = 3.\]

\begin{thm}
    The space of holomorphic \(1\)-forms on \(\GalCl\), \(\GalCl/\Subgr\) and \(\GalCl/\ShortKer\) is isomorphic to \(\Gamma \oplus \Gamma\), where \(\Gamma\) is the standard representation of \(S_3\), hence 
    \[q(\GalCl) = q(\GalCl/\Subgr) = q(\GalCl/\ShortKer) = 4\]
    and
    \[p_g(\GalCl) = 30, \quad p_g(\GalCl/\Subgr) = 12 ,\quad p_g(\GalCl/\ShortKer) = 6.\]
\end{thm}

\begin{proof}
    The computation of the irregularity follows the proof of Theorem 6.6 in \cite{BPS}. The geometric genus is then obtained as \(p_g = \chi(\cO) + q -1\).
\end{proof}

The computation of the invariants allows us to deduce that the degree $3$ maps $A/G\dashrightarrow \mathbb P^2$ and $A/K\dashrightarrow \mathbb P^2$ are different from the previously known degree $3$ map from abelian surfaces of types $(1,2)$ and $(1,6)$ to $\mathbb P^2$.
\begin{rem}
\label{irrationality}
    For a very general abelian surface $S$ of type $(1,2)$ or $(1,6)$ there are at least two different degree $3$ rational maps $S \dashrightarrow \mathbb P^2$ (non birationally equivalent). For an abelian surface $S$ of type $(1,2)$ the first degree $3$ map to $\mathbb P^2$ was constructed in \cite{yoshi}. The map we construct for $S=A/G$ above is different from the one constructed in \cite{yoshi} since the invariants of the Galois closure are different from the ones computed in \cite{BPS}. The same reasoning applies for the $(1,6)$ abelian surface.

    Let us remark that by \cite{martin} the polarization of a very general abelian surface having a degree $3$ morphism to $\mathbb P^2$ must be $(1,d)$ with $d=1,2,3,6$.
\end{rem}
\begin{rem}
 {The surface $X=W/\ShortKer$ is not the Galois closure of the map $\tilde \varphi/K: \tilde A/\ShortKer\to \mathbb P^2/\ShortKer.$} We will show in \autoref{prop:AlbaneseMap} that $X\to (A\times A)/K$ is the Albanese map of $X$. In particular, the Galois closure will be the image of $X$ in $\tilde A/K\times \tilde A/K$ under the degree $9$ isogeny $(\tilde A\times \tilde A)/K\to \tilde A/K \times \tilde  A/K $ which has many singularities. These singularities lie above the singular points of $\mathbb P^2/K$. In this case the Galois closure is singular because the codomain of $\tilde \varphi '$ is singular.
\end{rem}
\section{The geometry of the Albanese map} \label{sec:AlbaneseMap}
In this section we show that $X \to (A \times A)/\ShortKer$ is the Albanese map of \(X\) and we explore its local and global geometry. We shall assume that $A$ is very general.


First, let us observe the following.

\begin{lem}
    The self-intersection of \(X\) in \((A\times A)/\ShortKer\) is \(X^2 = 14\). 
\end{lem}
\begin{proof}
The cohomology class of \(W\) in \(A \times A\) can be easily computed (see for instance the proof of Theorem 1.1 in \cite{martin}). It is
    \[[W] = 3([A \times 0_A]+[0_A \times A]) - [\Delta_A] + S_{L_A}.\]
    Here \(S_{L_A} \coloneqq \pi_1^*(c_1(L_A))\cdot \pi_2^*(c_1(L_A))\) with $\pi_i:A \times A \to A$ denote the two projections.
    By computing the self-intersection of this class one gets 
    \[
    \begin{split}
    W^2 &= 18 ([A \times 0_A] \cdot [0_A \times A]) - 12 ([A \times 0_A] \cdot \Delta_A) - 2 (S_{L_A} \cdot \Delta_A) + (S_{L_A})^2 =\\
    &= 18 - 12 - 2 \cdot 12 + 144 = 126    
    \end{split}
    \]
    Now, \(W \subset A \times A\) is invariant with respect to the action of the diagonal part of \(K \times K < A \times A\) and its image via the quotient \(A \times A \to (A \times A)/K\) is \(X \subset (A \times A)/K\) ), so   
    \[X^2 = \frac{W^2}{9} = 14.\]
\end{proof}

\begin{prop} \label{prop:AlbaneseMap}
    The map $\alpha:X\to (A\times A)/K$ is the Albanese map of $X$. Moreover $\alpha$ is a local immersion and there are $2$ unique distinct points $q_1,q_2\in X$ such that $\alpha(q_1)=\alpha(q_2)$. 
\end{prop}

\begin{proof}
Let us check the second part of the Proposition. 
By \autoref{rem:Wsmooth} $W$ is embedded in $\tilde{A}\times \tilde{A}$. We study the restriction of the map \(\pi \colon \tilde{A} \times \tilde{A} \to A \times A\) to \(W\). We note that, for any exceptional divisor \(E\) in \(\tilde{A}\), \(\pi\) contracts \(\{a\} \times E\) to \((a,p_E)\) and \(E \times \{a\}\) to \((p_E,a)\) for any \(a \in A \smallsetminus \{p_E\}\), and maps \(E \times E\) to \((p_E,p_E)\). If two points \((x_1, x_2)\) and \((y_1, y_2)\) of \(\tilde{A} \times \tilde{A}\) have the same image via \(\pi\), then either\\
\begin{minipage}{0.4 \textwidth}
    \begin{enumerate}
        \item \(x_i = y_i \eqqcolon a\) and \(x_j, y_j \in E\)
    \end{enumerate}
\end{minipage}
or 
\begin{minipage}{0.4 \textwidth}
    \begin{enumerate}
    \setcounter{enumi}{1}
        \item \(x_1, x_2, y_1, y_2 \in E\). 
    \end{enumerate}
\end{minipage}\\
We, now, recall from \autoref{GenusHyperellCurves} and \autoref{QuotientRamification} that, for any \(z \in E\), \({\tilde{\varphi}}^{-1} (\tilde{\varphi}(z)) \smallsetminus \{z\} \subset H_E \smallsetminus E\), and \({\tilde{\varphi}}^{-1} (\tilde{\varphi}(z)) = \{z\}\) if and only if \(z \in E \cap H_E\).\\
Therefore, if \((x_1, x_2), (y_1, y_2) \in W \subset \tilde{A} \times \tilde{A}\), then, in case (1), \(\tilde{\varphi}(x_j) = \tilde{\varphi}(a) = \tilde{\varphi}(y_j)\), which means that \(x_j = y_j\), while, in case (2), \(\tilde{\varphi}(x_1) = \tilde{\varphi}(x_2)\) and \(\tilde{\varphi}(y_1) = \tilde{\varphi}(y_2)\), so \(x_1 = x_2 \eqqcolon x_E \in E \cap H_E\) and \(y_1 = y_2 \eqqcolon y_E \in E \cap H_E\). 


Since \(E \cap H_E\) consists in two points \(\{x_E,y_E\}\), the points \((x_E, x_E)\) and \((y_E, y_E)\) are the only ones contracted by \(W \to A \times A\).

Let us denote with $v_{x_E},v_{y_E}$ the tangent directions to the node of the hyperelliptic curve in $A$. In order to finish the proof of the second part of the statement it is enough to observe that $d\alpha_{(x_E,x_E)}T_{(x_E,x_E)}X=\langle (v_{x_E},0),(0,v_{x_E})\rangle$ and $d\alpha_{(y_E,y_E)}T_{(y_E,y_E)}X=\langle (v_{y_E},0),(0,v_{y_E})\rangle$. So as long as the hyperelliptic curve is nodal in $A$ (as it happens for very general $A$) we get that $\alpha$ is a local immersion and moreover \(d\alpha_{(x_E,x_E)}T_{(x_E,x_E)}X\oplus d\alpha_{(y_E,y_E)}T_{(y_E,y_E)}X =  T_{(p_E,p_E)}(A \times A)\).

Taking the quotient by \(\ShortKer\) all the \((x_E, x_E)\) and \((y_E, y_E)\) are mapped respectively to the points \(q_1\) and \(q_2\). Since \(W\) is embedded in \(\tilde{A} \times \tilde{A}\), \(X\) is embedded in \((\tilde{A} \times \tilde{A})/\ShortKer\) and \(q_1\) and \(q_2\) are the only points of \(X\) such that \(\alpha(q_1) = \alpha(q_2)\). In particular, $\alpha(X)$ has a unique singular double point.

We, now, check the first part of the statement. Let \(X\to \mathrm{Alb}(X)\) be the Albanese map of \(X\). Suppose, by contradiction, that there is an isogeny $\psi:\mathrm{Alb}(X) \to (A\times A)/\ShortKer$ of degree $d \ge 2$. On the one hand, the preimage through \(\psi\) of the singular locus of \(X \subset (A \times A)/\ShortKer\) has degree \(d\). On the other hand, it contains the self intersection of $\cup_{g\in G}g^*X$ where $\mathrm{Alb}(X)/G=(A\times A)/\ShortKer$ for a group of translations $G$ of cardinality $d$. We get  $(\cup_{g\in G}g^*X)^2=\binom{d}{2}X^2>d$, where the latter self intersection is in \(\mathrm{Alb}(X)\). This gives a contradiction.
\end{proof}

As a corollary, we can compute the Albanese map of $W$.

\begin{cor}
    The map $W\to A\times A$ is the Albanese map of $W$.
\end{cor}

\begin{rem}
\label{hodge}
    In \cite{schoen2007}, Schoen constructs a singular surface \(V\) with low invariants and proves that its Albanese variety \(A\) admits a deformation for which the cohomology class \([V]\) remains of Hodge type. The first step to prove the result is to show that the semi-regularity map $H^1(V, \mathcal N_{V/A}) \to H^3(\Omega^1_A,A)$, defined in \cite{bloch}, is injective. In the paper of Schoen this follows from the diagram \[
    \begin{tikzcd}
      H^1(T_A,A) \arrow{dr}\arrow{r} & H^1(V, \mathcal N_{V/A}) \arrow{d} \\
       & H^3(\Omega^1_A,A).
    \end{tikzcd}\]
    
   where the diagonal map is given by the cup product by \([V]\). In fact, $h^1(V, \mathcal N_{V/A})$ and the rank of the cup product map are both equal to $14$. This idea might apply to \(\alpha \colon X \to \mathrm{Alb}(X) = (A\times A)/K\) studied in \autoref{prop:AlbaneseMap}, for which we can show \(h^0(\mathcal{N}_{\alpha}) = h^2(\mathcal{N}_{\alpha}) = 4\) and \(h^1(\mathcal{N}_{\alpha}) = 14\).  However, Bloch's theory holds for local complete intersections, and the singularity of \(\alpha(X)\) in \(\mathrm{Alb}(X)\) is locally the intersection of two planes in \(\bC^4\), which is not a local complete intersection; hence, showing the injectivity of the semi-regularity map will require some work into a more general setting.  In this generality $h^1$ of the normal bundle should be replaced by $\mathrm{Ext}^2_{(A\times A)/K}(\mathcal O_{\alpha(X)},\mathcal{O}_{\alpha(X)})$. This was proven by Buchweitz and Flenner in \cite{semiregolarita2}. In our case, the rank of the cup product can be computed as $13$, so $\alpha(X)$ would be semi-regular in its Albanese if $\mathrm{ext}^2_{(A\times A)/K}(\mathcal O_{\alpha(X)},\mathcal{O}_{\alpha(X)})=13$. Once semiregularity is established, it would remain to determine the locus of abelian varieties where the cohomology class of $X$ stays of Hodge type (it is not clear to us whether this locus is of larger dimension of the locus of abelian four-folds isogenous to a product).
\end{rem}

\section{The surface \texorpdfstring{$X$}{X} admits no fibration to curves of positive genus}
\label{Section: 7}
In this section, we show that the surfaces $W$ and $X$ do not admit fibrations to curves of positive genus. This implies, in particular, that they are not quotients of products of curves.
\begin{thm}
    Suppose $\mathrm{NS}(A)=\mathbb Z \cdot c_1(L_A)$, then the surface $X$ does not map onto any curve of positive genus.
\end{thm}
\begin{proof}
    Suppose we have a morphism $\psi \colon X \to C$. We can suppose $C$ of genus $\ge 2$ since $A$ is simple and $\mathrm{Alb}(X)$ is isogenous to $A \times A$. By Castelnuovo-de Franchis Theorem having a fibration onto a curve of genus $\ge 2$ is equivalent to have a simple tensor $\eta$ in the kernel of
    \[
    \lambda \colon \bigwedge ^2 H^0(\Omega^1_X)\to H^0(\omega_X).
    \] 
    The wedge product morphism is equivariant with respect to the action of the Galois group $S_3$. Recall that the representation of $S_3$ on $H^0(\Omega^1_X)$ is $\Gamma \otimes H^0(\Omega ^1_{A/K})$, where $\Gamma$ is the standard representation. (cf. \cite{BPS}, Theorem 6.6). Now the orbit $\{g\cdot \eta \}$ in $\mathbb{P}\bigwedge ^2 H^0(\Omega^1_{A/K})$ is of cardinality either $1$, $2$, $3$ or $6$. Let us divide the proof according to the cardinality of the orbit.
    \begin{itemize}
        \item Suppose the cardinality of the orbit is one. We get that $\eta\in \bigwedge ^2 \langle \omega \rangle \otimes \Gamma$, for some $\omega \in H^0(\Omega^1_{A/K})$. This means that $\psi^{*} H^0(\Omega^1_C) \cap H^0(pr_1^{*} \Omega^1_{A/K})$ is one dimensional. Hence, it would induce a non-trivial sub-hodge structure on $H^0(\Omega^1_{A/K})$ contradicting the simplicity of $A/K$. 
        \item Suppose the cardinality of the orbit is 2. We get two simple tensors in the kernel. If the line joining them is a line of simple tensors  we get a map $X\to C$ for a curve of genus $ 3,4$ (cannot be higher for irregularity reasons). If the genus of the curve is $3$ we get that $(A\times A)/K\to \mathrm{Jac}(C)$ has an elliptic curve as kernel, contradicting the simplicity of $A$. If $C$ is of genus $4$ it is easy to see that the Albanese dimension of $X$ would be $1$, contradiction. Hence, we may suppose that the line joining the points in the orbit of $\eta$  is meeting transversally the Grassmannian of simple tensors. Let $\sigma :X \to X$ be an involution in $S_3$, we have that $\sigma ^*\eta$ is the other point in the orbit of $\eta.$ This means that we get a dominant map $X\to C\times C$, where the second projection is obtained by the original fibration pre-composed with $\sigma$. This implies that $C$ is of genus $2$. This situation leads to the following commutative diagram:
        \[
        \begin{tikzcd}
            X \arrow[r, "\sigma"] \arrow[d] & X \arrow[d]\\
            C \times C \arrow[r, "\sigma'"] & C \times C.
        \end{tikzcd}
        \]
        This induces a map $X/\sigma \simeq \tilde{A}/K\to (C\times C)/\sigma'\simeq \mathrm{Sym}^2C$, which induces an isogeny $A/K \to \mathrm{Jac}(C)$ (of degree $\ge 2$, since $A/K$ is a very general $(1,6)$ abelian surface and $\mathrm{Jac}(C)$ is principally polarized). By the universal property of the fibered product we get a birational map \[
        \tilde A/K \to \mathrm{Sym}^2(C)\times_{\mathrm{Jac}(C)}A /K.\]
        Since $\mathrm{Sym}^2C\to \mathrm{Jac}(C)$ contracts a $(-1)$ curve and $A/K\to \mathrm{Jac}(C)$ is étale of degree $\ge 2$ we deduce that $\mathrm{Sym}^2(C)\times_{\mathrm{Jac}(C)} A/K$ contains at least two $(-1)$ curves. Hence $\tilde{A}/K$ contains at least two $(-1)$ curves, contradiction. 
        \item Suppose the cardinality of the orbit is $3$ or $6$. Let us consider the linear  space $V\subset \mathrm{ker}(\lambda)$ spanned by the orbit. We claim that $\bP \mathrm{ker}(\lambda )\cap \mathrm{Gr}(2,4)$ is at least $1-$dimensional\footnote{Here $\mathrm{Gr}(2,4)$ denotes the Grassmanian of simple tensors.}. Indeed, if $\bP V$ is a line, then $\bP V\subset \mathrm{Gr}(2,4)$, since the Grassmanian is a quadric. Now $\bP \mathrm{ker}(\lambda)$ cannot contain any line, otherwise we would get a map $X\to C$, with $C$ a curve of genus $3,4$ (this is impossible, see the above bullet point). If $\bP \mathrm{ker}(\lambda)\cap \mathrm{Gr}(2,4) $ is $1-$dimensional and does not contain any line, then we get a positive dimensional family of maps $X\to C_t$, where $C_t$ is a smooth genus $2-$curve. This leads to a contradiction because the morphism of Hodge structures $H^0(\Omega^1_{C_t})\to H^0(\Omega^1_X)$ is rigid.
    \end{itemize}
\end{proof}
 \begin{rem}
 Having a map onto a curve of genus $\ge 2$ is a topological invariant. So if $A$ is isogenous to a product of elliptic curves it is still true that $X$ is not fibered over curves of genus $\ge 2$ (of course it is fibered over elliptic curves in this case).
 \end{rem}
This has some standard consequences on the wedge product of $1-$forms on the surface $X$.
\begin{cor} \label{cor:RankLambda}
    The map \[
\bigwedge^2 H^0(\Omega^1_X,X)\to H^0(\omega_X)
\] has rank $5$ and has one non-simple tensor in the kernel.
\end{cor}
\begin{proof}
    By Castelnuovo-de Franchis theorem the map above has a non-simple tensor in the kernel if and only if $X$ is fibered over a curve of genus $\ge 2$. Hence the wedge product does not have any simple tensor in the kernel. This implies that the rank is at least $5$ (the kernel cannot intersect $\mathrm{Gr}(2,4)$).

    It is easy to observe that the kernel is non-trivial. In fact, the (birational) quotient of $(W,A \times A)$ by the action of the Galois group $S_3$ is $(\mathrm{Kum}_2(A),\mathbb P^2)$. Hence, the $2$-form $\eta$ pulled back from $\mathrm{Kum}_2(A)$ has to vanish on $W$. This $2$-form is invariant under the diagonal action of $K$, hence we have a non-trivial kernel of the wedge product of $1-$forms on $X$, as desired. 
\end{proof}
We can also deduce that the curve $X$ is not a quotient of a product of curves.
\begin{cor}
    The surface $X$ is not birational to $(C_1\times C_2)/G$.
\end{cor}
\begin{proof}
The irregularity of \((C_1 \times C_2) / G\) is the sum of the genus of \(C_1/G\) and \(C_2/G\), but this means that the genus of either \(C_1/G\) or \(C_2/G\) is \(\ge 2\), which is impossible by \autoref{cor:RankLambda}.
\end{proof}
\section{The canonical image of \texorpdfstring{$X \coloneqq \GalCl/\ShortKer$}{the quotient by K}} \label{sec:CanonicalImage}
A natural question is the nature of the canonical map of the surfaces we are studying. This section aims to show that the canonical map of \(X \coloneqq \GalCl/\ShortKer\) has degree 2 on a degree \(12\) surface in \(\bP^5\) with \(44\) nodes.

First, we need to find the intersection of \(X\) with the fixed locus of the map \((-1)_{(\tilde{A} \times \tilde{A})/\ShortKer}\) induced by $(-1)_{\tilde{A} \times \tilde{A}}$. The ramification of the degree \(3\) map \(\tilde{\varphi} \colon \tilde{A} \to \bP^2\) contains fourteen fixed points of \((-1)_{\tilde{A}}\). We call them \(\{p_i\}_{i = 1, \dots, 14}\) two of which, say \(p_{13}\) and \(p_{14}\), have multiplicity \(3\). For \(i = 1, \dots, 12\), \((\tilde{\varphi})^{-1}(\tilde{\varphi}(p_i)) = \{p_i, q_i\}\) for some \(q_i\) in the exceptional divisor. Finally, there are three other fixed points of \((-1)_{\tilde{A}}\), we call them \(\{t_1, t_2, t_3\}\), which are sent to the same point by \(\tilde{\varphi}\). Hence, \(X\) contains the classes of the following \(44\) fixed points of \((-1)_{\tilde{A} \times \tilde{A}}\), which are distinct fixed points of \((-1)_{(\tilde{A} \times \tilde{A})/\ShortKer}\) because the group \(\ShortKer\) has odd order:
\[W[2] \coloneqq \{(t_i, t_j) | i \neq j\} \cup \{(p_i, p_i) | i = 1, \dots, 14\} \cup \{(p_i, q_i), (q_i, p_i) | i = 1, \dots, 12\}.\]

In what follows we call:
\begin{itemize}
    \item \(X[2]\) the image of \(W[2]\) via the projection \(\tilde{A} \times \tilde{A} \to (\tilde{A} \times \tilde{A})/\ShortKer\);
    \item \(X_{44}\) the blow up of \(X\) along \(X[2]\);
    \item \(\eta \colon X_{44} \to X\) the blow up;
    \item \(T \coloneqq X_{44} / (-1)_{(\tilde{A} \times \tilde{A})/\ShortKer}\);
    \item \(\pi \colon X_{44} \to T\) the projection.
\end{itemize}



\begin{thm}\label{thm:CanonicalMapX}
    The canonical map of \(X \coloneqq \GalCl/\ShortKer\) has degree $2$ onto its image and its image has \(44\) nodes.
\end{thm}

\autoref{thm:CanonicalMapX} follows from \autoref{prop:CanMapXFactorizesT} and \autoref{prop:CanMapTBirational}.

\begin{prop}\label{prop:CanMapXFactorizesT}
	The canonical map of \(X\) factorizes through the canonical map of \(T\).
\end{prop}

\begin{proof}
    To prove the first part of the statement, it suffices to show that the canonical map of \(X\) factorizes through that of \(T\). This is true if \(p_g(T) = p_g(X) = 6\). \\
    Notice that, applying the Hurwitz formula to \(\eta\) and \(\pi\), we get two expressions of the Canonical divisor of \(X_{44}\):
    \[K_{X_{44}} \sim \eta^{*} K_X + \sum_{i=1}^{44} F_i\]
    \[K_{X_{44}} \sim \pi^{*} K_T + \sum_{i=1}^{44} F_i\]
    where the \(F_i\)'s are the exceptional divisors of \(X_{44}\).\\
    This implies that \(\pi^{*} K_T = \eta^{*} K_X\), which allows us to compute the self-intersection of \(K_T\):
    \[K_T^2 = \frac{1}{2} (\pi^{*} K_T)^2 = \frac{1}{2} (\eta^{*} K_X)^2 = \frac{1}{2} K_X^2 = 12.\]
    Now, by Persson's formulas on double covers,
    \[c_2(X_{44}) = 2 c_2(T) + 2 p_a\left(\sum_{i =1}^{44} F_i\right) - 44 \cdot 2,\]
    which, since \(c_2(X_{44}) = c_2(X) + 44 = 56\) and the exceptional divisors are rational, implies that \(c_2(T) = 72\).\\
    Then
    \[\chi(\cO_T) = \frac{1}{12}\left(c_1(T)^2 + c_2(T)\right) = 7;\]
    and, finally, \(\chi(\cO_T) = 1 + q(T) + p_g(T) = 1 + 0 + p_g(T)\).
\end{proof}

We now explain a bit of the geometry of the canonical map of \(X\) and \(T\). In order to do so we need some properties of the surface \(Z \coloneqq X/A_3\) and, more specifically, of \(S \coloneqq \tilde{Z}/(-1) \sim_{\mbox{bir}} T/A_3\), which has been studied by the second author in \cite[Theorem 6.13, Theorem 6.14]{grossi2025galois}.

\begin{prop}
    \(S\) is an Horikawa surface of invariants:
    \[%
	c_2(S) = 56, \quad%
	K_{S}^2 = 4, \quad%
	q(S) = 0, \quad%
	p_g(S) = 4%
	\]
    In particular, \(S\) has canonical map of degree \(2\) onto a quadric in \(\bP^3\).
\end{prop}

\begin{prop}\label{prop:CanMapTBirational}
	The canonical map of \(T\) is birational onto its image and contracts the \(44\) \(-2\)-curves.
\end{prop}

\begin{proof}
We will consider the following diagram.
\[
\begin{tikzcd}
    T \arrow{r}{\varphi_{K_T}} \arrow{d}{\psi} & \mathbb P^5 \arrow[dashed,d] \\
    S \arrow{r}{\varphi_{K_S}} & \mathbb P^3
\end{tikzcd}
\]
We call \(R_S\) the ramification of \(\varphi_{K_S}\), \(\sigma\) the involution on \(S\) inducing \(\varphi_{K_S}\) and \(\tau\) the order \(3\) automorphism inducing \(\psi\).

    We divide the proof in 4 steps.
    \begin{enumerate}
        \item \emph{The canonical map of \(T\) does not factorize through \(S\).} By contradiction, suppose that the canonical map of \(T\) factorizes through \(S\) and let's call \(L\) the line bundle defining the map \(S \to \bP^5\). Then, \(\psi^* |H^0(L)|\) are the canonical divisors of \(T\) invariant under the action of \(A_3\). Since \(h^0(\omega_T)^{A_3} = h^0(\omega_S) = p_g(S) = 4 < 6 = h^0(\omega_T)\), the action of \(A_3\) on \(|H^0(\omega_T)|\) is not trivial, which gives a contradiction. 
        \item \emph{\(\varphi_{K_T}\) has degree either \(1\) or \(2\).} Since the canonical map of \(T\) does not factorize through \(S\), there is a map of degree at least \(3\) between \(\varphi_{K_T}(T)\) and \(\varphi_{K_S}(S)\). Hence, the degree of \(\varphi_{K_T}(T)\) is \(\geq 3 \cdot \deg(\varphi_{K_S}(S)) = 6\). Since \(K_T^2 = 12\) the claim follows.

        \item \emph{\(\varphi_{K_T}\) has degree \(1\).} By contradiction, suppose that \(\varphi_{K_T}\) has degree \(2\). Then, by the commutativity of the diagram, \(\varphi_{K_T}\) is induced by a lifting \(\tilde{\sigma}\) of \(\sigma\). Note that \(\langle\tilde{\sigma}, \tau\rangle \cong S_3 \subset \mbox{Aut}(T)\) as we will show in (4). Let us denote by \(g(\tau)\) and \(g(\tilde{\sigma})\) the representative of \(\tau\) and \(\tilde{\sigma}\) in \({\bP}GL(H^0(\omega_T))\). Since the canonical map \(T \to \varphi_{K_T}(T)\) is induced by \(\tilde{\sigma}\), we deduce that \(g(\tilde{\sigma})\) is the identity. On the other hand, since, as seen in point (1), the action of \(A_3\) is nontrivial, we get that \(g(\tau)\) is a nontrivial element of order \(3\). We get a contradiction because \(\tilde{\sigma} \circ \tau\) is an element of order \(2\) while \(g(\tilde{\sigma} \circ \tau) = g(\tilde{\sigma}) g(\tau) = g(\tau)\) has order \(3\).

        \item {\(\langle\tilde{\sigma}, \tau\rangle \cong S_3 \subset \mbox{Aut}(T)\).} Since \(\langle\tilde{\sigma}, \tau\rangle\) has order \(6\), it suffices to show that \(\tilde{\sigma}\) and \(\tau\) do not commute. This, in turn, is implied by the existence of \(p \in \mathrm{Fix}(\sigma) \subset S\) such that \(\tilde{\sigma}\) acts nontrivially on \(\psi^{-1}(p)\). Let \(C\) be a smooth canonical section of \(S\) avoiding the image of the fixed points of \(A_3\) in \(T\). The image of \(C\) is a quadric, therefore \(\sigma\) acts as an hyperelliptic involution on \(C\). Since the genus of \(C\) is \(5\), it has \(12\) Weierstrass points. The pullback of \(C\) via \(\psi\) is a smooth curve \(\tilde{C}\) of genus \(13\) containing \(V_{36}\) the set of the \(36\) preimages of the Weierstrass points. We get an induced action of \(\tilde{\sigma}\) on this curve and on \(V_{36}\). This action cannot fix all the points of \(V_{36}\) by the Hurwitz formula applied to \(\tilde{C} \to \tilde{C}/\tilde{\sigma}\). We, then, can choose \(p \in \mathrm{Fix}(\sigma)\) as the image of a point in \(V_{36}\) not fixed by \(\tilde{\sigma}\).
    \end{enumerate}
    
\end{proof}
    
\bibliography{refer}
\bibliographystyle{myalphaspecial}
    
\end{document}